\documentclass[10pt, twocolumns]{IEEEtran}

\usepackage{cite}
\usepackage{amsmath,amssymb,amsfonts}
\usepackage{algorithmic}
\usepackage{textcomp}
\usepackage{color}
\usepackage{bm}
\usepackage{graphicx}
\usepackage{enumerate}
\usepackage{multicol}
\usepackage{booktabs,tabularx,makecell,array}

\newcolumntype{L}{>{\raggedright\arraybackslash}X}

\newtheorem{lemma}{Lemma}
\newtheorem{assumption}{Assumption}

\newtheorem{remark}{Remark}

\newtheorem{theorem}{Theorem}

\newtheorem{problem}{Problem}

\title{
\huge{
Dynamic Nash Equilibrium Seeking for a Class of Nonlinear Uncertain  Multi-agent Systems
}
}

\author{Weijian~Li, and Yutao~Tang

\thanks{W.~Li is with the School of Civil and Environmental Engineering, Cornell University, Ithaca, NY,
14850, USA. E-mail: 
\texttt{wl779@cornell.edu}.
}

\thanks{Y.~Tang is with  the School of Intelligent Engineering and	Automation, Beijing University of Posts and Telecommunications, Beijing
100876, China. E-mail: 
\texttt{yttang@bupt.edu.cn}.
}
}

\begin{document}

\maketitle

\begin{abstract}                        
We consider seeking a Nash equilibrium (NE) of a monotone game, played by dynamic agents which are modeled as a class of lower-triangular nonlinear uncertain dynamics with external disturbances.
We establish a general framework that converts the problem into a distributed robust stabilization problem of an appropriately augmented system.
To be specific, we construct a virtual single-integrator multi-agent system, as a reference signal generator, to compute an NE in a fully distributed manner.
By introducing internal models to tackle the disturbances, as well as embedding the virtual system,
we derive an augmented system.
Following that, we show that the outputs of all agents reach an NE of the game if the augmented system can be stabilized by a control law.
Finally, resorting to a backstepping procedure, we design a distributed state-feedback controller to stabilize the augmented system semi-globally.
\end{abstract}


\section{INTRODUCTION}

Distributed NE seeking for monotone games has received a flurry of research interest, motivated by its broad applications from network congestion control,  communication networks, smart grids to social networks \cite{saad2012game, ghaderi2014opinion, yin2011nash}.
The basic setup is that in a multi-agent system, each player (agent) aims to minimize a local cost function depending on its own strategy as well as on the strategies of its opponents. All players try to reach an NE, whereby no player can decrease its local cost by unilaterally changing its own decision.
A variety of distributed algorithms have been proposed over the years, including best-response,  gradient-play, payoff-based learning and operator splitting approaches \cite{hu2022distributed, ye2023distributed, belgioioso2022distributed}.
One of the most studied methods is the gradient-play scheme, which is easy to be implemented under full- and partial-decision information settings \cite{shamma2005dynamic, gadjov2018passivity, tatarenko2020geometric}.

In practice, agents may have inherent dynamics, and their strategies are outputs of a dynamical system.
Examples can be found in coordination of mobile sensor networks \cite{stankovic2011distributed},
load allocation for plug-in electric vehicles \cite{ito2012disturbance}, and distributed control of wind farms \cite{marden2013model, krilavsevic2021learning}.
The agent dynamics have a great influence on the decision-making process, and thus, one should take them into consideration when developing distributed algorithms.
Recent research efforts have focused on this area.
For aggregative games, a proportional-integral feedback algorithm was explored for a class of second-order passive systems in \cite{de2019feedback}, 
and distributed gradient-based  protocols were introduced for 
Euler–Lagrange systems and nonlinear systems with unit relative degree in  \cite{deng2019distributed, zhang2019distributed}.
Monotone games played by dynamic agents were considered in \cite{romano2019dynamic}, and distributed NE seeking strategies with disturbance rejection were proposed.
The results were extended to distributed generalized NE computation of monotone games with convex separable coupling constraints in \cite{bianchi2021continuous}. In \cite{ye2020distributed}, control schemes with bounded inputs were further investigated.
However,  in \cite{romano2019dynamic, bianchi2021continuous, ye2020distributed}, the agent dynamics took special forms of multi-integrators.
For heterogeneous linear multi-agent systems, output feedback strategies were provided for quadratic games in \cite{guo2021linear}, resorting to linear output regulation.
In \cite{tang2022nash}, a class of high-order nonlinear systems with unknown dynamics was considered, and a distributed adaptive protocol was developed. Under switching topologies,
the distributed NE seeking problem was investigated in \cite{huang2024distributed} for a class of nonlinear systems with bounded disturbances.
By adaptive backstepping approaches, distributed NE seeking for a class of nonlinear uncertain systems was addressed in \cite{meng2024distributed}.
Besides, a multi-cluster game problem with agents modeled by second-order dynamics was explored in \cite{deng2024distributed}.
However, existing literature has not reported distributed NE seeking strategies for complex nonlinear multi-agent systems with both uncertainties and external disturbances.

Inspired by the above observations, we focus on designing a distributed control protocol to steer the outputs of a multi-agent system to an NE of a monotone game.
Our main contributions are summarized as follows.
First, we consider distributed NE seeking for a class of nonlinear multi-agent systems in a lower-triangular form, allowing both uncertain parameters and external disturbances.
The system covers those in \cite{romano2019dynamic, guo2021linear, zhang2019distributed, zhang2023distributed, deng2024distributed, huang2024distributed}
as special cases, and is discussed for the first time to the best of our knowledge.
Second, by constructing a virtual reference signal generator for NE computation and introducing internal models to handle disturbances, we establish a general framework that reformulates the problem as stabilizing an appropriately augmented system. 
Compared with the distributed design in \cite{romano2019dynamic, zhang2019distributed, deng2024distributed, meng2024distributed}, our method is more flexible since the NE seeking and reference tracking problems are solved separately. 
In contrast to \cite{tang2020optimal, zhang2023distributed, tang2022nash}, we indicate that the framework can solve NE seeking problems for complex nonlinear systems with uncertainties and disturbances.
Last but not least,  by backstepping techniques, we show that a linear distributed state-feedback controller can be employed to solve the problem. Distinct from \cite{huang2004general, su2014cooperative, su2016cooperative}, our method tackles reference tracking, as well as NE seeking.

This paper is organized as follows. In Section \ref{sec:formulation}, we introduce necessary preliminaries, and formulate the problem. Then we establish a general framework in Section \ref{sec:conversion}, and present our main results in Section \ref{sec:results}. 
In Section \ref{sec:example}, we provide an illustrative example.
Finally, we give concluding remarks in Section \ref{sec:conclusion}.

\section{PRELIMINARY AND FORMULATION}
\label{sec:formulation}

In this section, we introduce some necessary concepts and formulate the distributed NE seeking problem.

\subsection{Mathematical Preliminary}

Let $0_m$ ($1_m$) be the $m$-dimensional column vector with all entries of $0$ ($1$), and $I_n$ be the $n$-by-$n$ identity matrix.
We simply write $\mathbf{0}$ for vectors of zeros with appropriate dimensions when there is no confusion.
Let $(\cdot)^\top$, $\otimes$ and $\Vert \!\cdot\! \Vert$ be the transpose, the Kronecker product and the Euclidean norm.
Let $X \times Y$ be the Cartesian product of sets $X$ and $Y$. Given $x_i \in \mathbb R^{n_i}$,
${\rm col}\{x_1, \dots, x_N\} = [x_1^{\top}, \dots, x_N^{\top}]^{\top}$.
The compact set $\bar {\mathbb Q}_R^s$ is defined as $\bar {\mathbb Q}_R^s = \big\{y = {\rm col}\{y_1, \dots, y_s\} \in \mathbb R^s: |y_j| \le R, j \in \{1, \dots, s\} \big\}$.
For a positive definite and radically unbounded function $V: \mathbb R^n \rightarrow \mathbb R$, 
the compact set $\bar \Omega_c(V(x))$ is defined as
$\bar \Omega_c(V(x)) = \{x \in \mathbb R^n: V(x) \le c \}$, and the open set $\Omega_c(V(x))$ is defined as  $\Omega_c(V(x)) = \{x \in \mathbb R^n: V(x) < c \}$. 

An operator 
$F: \mathbb R^n \rightarrow \mathbb R^n$ is monotone  if $\langle x - y, F(x) - F(y)\rangle \ge 0, \forall x, y \in \mathbb R^n$,  $\underline l$-strongly monotone if $\langle x - y, F(x) - F(y)\rangle \ge \underline l \Vert x - y \Vert^2, \forall x, y \in  \mathbb R^n$, and $\overline l$-Lipschitz continuous if $\Vert F(x) - F(y)\Vert \le \overline l \Vert x - y \Vert, \forall x, y \in  \mathbb R^n$.

Consider a multi-agent network modeled by an undirected graph $\mathcal G(\mathcal I, \mathcal E, \mathcal A)$, where $\mathcal I = \{1, \dots, N\}$ is the node set, $\mathcal E \subset \mathcal I \times \mathcal I$ is the edge set, and $\mathcal A = [a_{ij}] \in \mathbb R^{N \times N}$ is the adjacency matrix such that $a_{ij} = a_{ji} > 0$ if $(i, j) \in \mathcal E$, and $a_{ij} = 0$ otherwise. 
The Laplacian matrix $\mathcal L$ is $\mathcal L = \mathcal D - \mathcal A$, where $\mathcal D = {\rm diag}\{d_i\}$, and $d_i = \sum_{j \in \mathcal I} a_{ij}$.
The graph $\mathcal G$ is connected if there exists a path between any pair of distinct nodes.

\subsection{Problem Statement}

Consider a nonlinear multi-agent system composed of $N$ agents. The dynamics of agent $i$ is  described by 
\begin{equation}
\begin{aligned}
\label{agent:dyn}
\dot z_i &= f_{0i}(z_i, x_{1i}, v, w) \\
\dot x_{1i} &= f_{1i}(z_i, x_{1i}, v, w) + x_{2i} \\
&\vdots \\
\dot x_{r i} &= f_{r i}(z_i, x_{1i}, \dots, x_{r i}, v, w) + u_i \\
y_i &= x_{1i}, ~i \in \mathcal I
\end{aligned}
\end{equation}
where $\mathcal I = \{1, \dots, N\}$,  $z_i \in \mathbb R^{n_{z_i}}$ and $x_i \triangleq {\rm col}\{x_{1i}, \dots, \\
x_{r i}\} \in \mathbb R^{r}$ are the states, $u_i \in \mathbb R$ is the control input, $y_i \in \mathbb R$ is the output,  $w \in \mathbb W$ represents the parameter uncertainty, $v \in \mathbb R^{n_v}$ is the disturbance generated by an exosystem as
\begin{equation}
\label{exosystem}
\dot v = S v, ~v(0) \in \mathbb V_0
\end{equation}
both $\mathbb W \subset \mathbb R^{n_w}$ and $\mathbb V_0 \subset \mathbb V_0 \in \mathbb R^{n_v}$ are compact,
and moreover, the functions $f_{0i}$ and $f_{si}, s \in \{1, \dots, r\}$ are sufficiently smooth with $f_{0i}(\mathbf{0}, 0, \mathbf{0}, w) = 0$ and $f_{si}(\mathbf{0}, 0, \dots, 0, \mathbf{0}, w) = 0$ for all $w \in \mathbb W$.

All agents play an $N$-player noncooperative game, denoted by $\mathbf{G}(\mathcal I, J_i, \mathbb R)$. 
Specifically, agent $i$ is endowed with a local cost function $J_i(y_i, y_{-i}): \mathbb R^N \rightarrow \mathbb R$, where $y_i \in \mathbb R$ is its output strategy specified by (\ref{agent:dyn}), and 
$y_{-i} = [y_1, \dots, y_{i-1}, y_{i+1}, \dots,$ $y_N] \in \mathbb{R}^{N-1}$
denotes the strategy profile of its opponents.
Each agent changes its output according to (\ref{agent:dyn}) by choosing its control input, and moreover, communicates with its neighbors through an undirected graph $\mathcal G(\mathcal I, \mathcal E, \mathcal A)$.
All agents try to reach a steady-state output profile,  defined as an NE of $\mathbf{G}$ in this paper.
Given $\mathbf{G}(\mathcal I, J_i, \mathbb R)$, the profile $y^*  =  {\rm col}\{y_1^*, \dots, y_N^*\}$ is an NE  if $y_i^*  \in {\rm argmin}_{y_i} J_i(y_i, y_{-i}^*), ~\forall i \in \mathcal I.$

The controller $u_i$ is expected to take the form of
\begin{equation}
\label{control}
\left\{
\begin{aligned}
\dot \varrho_i &= \Xi_{1i}\big(\nabla_i J_i, x_j, \varrho_j\big) \\
u_i &=  \Xi_{2i}\big(\nabla_i J_i, x_j, \varrho_j\big),  
~j \in \mathcal I_i \cup \{i\}
\end{aligned}
\right.
\end{equation}
where $\varrho_i \in \mathbb R^{n_{\varrho_i}}$, $\Xi_{1i}$ and $\Xi_{2i}$ are sufficiently smooth functions to be specified,
$\nabla_i J_i(y_i, y_{-i}) \!=\! \partial J_i(y_i, y_{-i})/ \partial y_i$,
and $\mathcal I_i$ is the neighbor set of agent $i$, i.e., $\mathcal I_i = \{j | (i, j) \in \mathcal E\}$.
Let $x_c = {\rm col}\{z_1, x_1, \varrho_1, \dots, z_N, x_N, 
\varrho_N\}$ and $n_c = \sum_{i \in \mathcal I}(n_{z_i} + r + n_{\varrho_i})$.
Then we formulate the problem as follows.
\begin{problem}
\label{prb}
Consider the multi-agent system (\ref{agent:dyn}) and the  exosystem (\ref{exosystem}) under the undirected graph $\mathcal G$ with local functions $J_i$. Given any real number $R > 0$ and  nonempty compact set $\mathbb W \times \mathbb V_0  \subset \mathbb R^{n_w + n_v}$ containing the origin,
determine a distributed controller in the form of (\ref{control}) such that for any ${\rm col}\{w, v(0)\} \in \mathbb W \times \mathbb V_0$ and $x_c(0)\in \bar{\mathbb Q}_R^{n_c}$,
\begin{enumerate}[a)]
\item the trajectory of the closed-loop system consisting of (\ref{agent:dyn}) and (\ref{control}) exists, and is bounded over $[0, \infty)$.

\item The agents' output satisfies
$\lim_{t \to \infty} y_i(t) = y_i^*, i \in \mathcal I$, where $y^* = {\rm col}\{y_1^*, \dots, y_N^*\}$ is an NE of $\mathbf{G}(\mathcal I, J_i, \mathbb R)$.
\end{enumerate}
\end{problem}

\begin{remark}
Distributed NE seeking for noncooperative games has been investigated in \cite{gadjov2018passivity, bianchi2021continuous, romano2019dynamic, ye2020distributed}, but (\ref{agent:dyn}) was restricted to be single or multiple integrators.
This paper considers the nonlinear multi-agent system (\ref{agent:dyn}) in a more general form, which covers linear systems, nonlinear systems with unity relative degree, etc \cite{guo2021linear, wang2015distributed, zhang2019distributed, zhang2023distributed, tang2020optimal, deng2024distributed, huang2024distributed}.
In practice, (\ref{agent:dyn}) appears in many benchmark systems, including Chua’s circuit, Lorenz system, Duffing equation, and Van del Pol oscillators.
Compared to \cite{tang2022nash,meng2024distributed}, we allow the presence of uncertainties and disturbances.
\end{remark}

\section{GENERAL FRAMEWORK}
\label{sec:conversion}

Construct a virtual multi-agent system as the abstraction of (\ref{agent:dyn}). Let 
all virtual agents play the game $\mathbf{G}(\mathcal I, J_i, \mathbb R)$, and dynamics of agent $i$ be 
\begin{equation}
\label{vir:sys}
\dot p_i(t) = \omega_i(t), ~ i \in \mathcal I
\end{equation}
where $\omega_i$ is the input, and $p_i$ is the output.
In fact, (\ref{vir:sys}) can be viewed as a reference signal generator to compute an NE. Suppose 
$p(t) = {\rm col}\{p_1(t), \dots, p_N(t)\}$ approaches an NE. Then Problem \ref{prb} can be solved by designing $u_i$ such that $y_i(t)$ track the trajectory of $p_i(t)$.

The idea motivates us to establish a framework that
converts Problem \ref{prb} into a distributed robust stabilization problem of an appropriately augmented system. The conversion consists of the following three steps.
First, construct the reference signal generator (\ref{vir:sys}) to seek an NE.
Second, design internal models to handle the disturbances generated by (\ref{exosystem}). The nonlinear system (\ref{agent:dyn}), the virtual system (\ref{vir:sys}) and the internal models together form an augmented system.
After a suitable coordinate transformation, the stabilizability of the augmented system implies the solvability of Problem \ref{prb}. Third,  
design a controller to stabilize the augmented system semi-globally.

\begin{remark}
The framework is motivated by the designs in \cite{ tang2020optimal, zhang2023distributed, tang2022nash}.
However, we extend the approaches from optimal output consensus to distributed NE seeking with complex agent dynamics characterized by  nonlinearity, uncertainties and disturbances.
Compared to the methods in \cite{romano2019dynamic, zhang2019distributed, deng2024distributed, meng2024distributed},
our framework is more flexible and reconfigurable since the distributed NE seeking problem and the trajectory tracking problem can be solved separately.
Different from the output regulation problems in \cite{huang2004general, su2014cooperative, su2015cooperative, su2016cooperative}, the main challenge is to let $y_i(t)$ track the trajectory $p_i(t)$ generated by a virtual system (\ref{vir:sys}) rather than an output trajectory generated by (\ref{exosystem}).
\end{remark}

\subsection{Reference Signal Generator}

Similar to \cite{gadjov2018passivity, romano2019dynamic, bianchi2021continuous}, we consider that the virtual system (\ref{vir:sys}) seeks an NE of $\mathbf{G} (\mathcal I, J_i, \mathbb R)$ in a fully distributed manner, i.e., agent $i$ only knows $J_i$, and receives data from its neighbors via $\mathcal G$.
In this case, the exact partial gradient $\nabla_i J_i(p_i, p_{-i})$ cannot be computed.
Let agent $i$ maintain a vector 
$\mathbf{p}_i = {\rm col}\{p^i_1, \dots, p^i_{i-1}, p_i, p^i_{i+1},\dots, 
p^i_N\}  \in \mathbb R^N$, where $p^i_k$ is agent $i$'s estimate of agent $k$'s strategy, and $p_i$ is its actual strategy. 
Define a pseudo-gradient mapping $F: \mathbb R^N \rightarrow \mathbb R^N$ as
$F(p) = {\rm col}\{\nabla_1 J_1(p_1, p_{-1}), \dots, \nabla_N J_N(p_N, p_{-N})\}$, and an extended pseudo-gradient mapping $\mathbf{F}: \mathbb R^{N^2} \rightarrow \mathbb R^N$ as $\mathbf{F}(\mathbf{p}) = {\rm col}\{\nabla_1 J_1(\mathbf{p}_1), \dots, \nabla_N J_N(\mathbf{p}_N)\}$, where $\mathbf{p} = {\rm col}\{\mathbf{p}_1, \dots, \mathbf{p}_N\}$, and
$\nabla_i J_i(\mathbf{p}_i) = {\partial J_i(\mathbf{p}_i)}/ \partial p_i$.
We further make the following assumptions \cite{gadjov2018passivity, tatarenko2020geometric, bianchi2021continuous}.

\begin{assumption}
\label{ass:convex}
For every $i \in \mathcal I$, $J_i$ is continuously differentiable and convex in $x^i$, given $x^{-i}$. 
Furthermore, both $F$ and $\mathbf{F}$ are $\overline l_F$-Lipschitz continuous, and $F$ is $\underline l_F$-strongly monotone.
\end{assumption}

\begin{assumption}
\label{ass:graph}
The undirected graph $\mathcal G$ is connected.
\end{assumption}

With these preparations, the following fully distributed gradient-play dynamics can be employed as a reference signal generator:
\begin{equation}
\label{gradient-play}
\left\{
\begin{aligned}
\dot p_i &= - \gamma_1\nabla_i J_i(\mathbf{p}_i) - \gamma_1 \gamma_2\sum\nolimits_{j \in \mathcal I_i} a_{ij}(p_i - p^j_i) \\
\dot p^i_{k} &= - \gamma_1\gamma_2\sum\nolimits_{j \in \mathcal I_i} a_{ij}(p^i_k - p^j_k),
~k \in \mathcal I \backslash \{i\}
\end{aligned}
\right.
\end{equation}
where $\gamma_1, \gamma_2 > 0$, and $a_{ij}$ is the $(i, j)$-th entry of the adjacency matrix of $\mathcal G$.
Let $\mathcal L$ be the Laplacian matrix of $\mathcal G$, and
$\mathbf{L} = \mathcal L \otimes I_N$. Define 
$\mathcal R_i = [0_{i -1}^{\top}, 1, 0_{n-i}^{\top}] \in \mathbb R^{1 \times N}$ and
$\mathcal R = {\rm blkdiag}\{\mathcal R_i\}_{i \in \mathcal I}$.
Then dynamics (\ref{gradient-play}) reads as
\begin{equation}
\label{gradient-play:sim}
\dot{\mathbf{p}} =  - \gamma_1 {\mathcal R}^{\top} \mathbf{F}(\mathbf{p}) - \gamma_1 \gamma_2 \mathbf{L} \mathbf{p}.
\end{equation}
The next lemma establishes the convergence of (\ref{gradient-play:sim}).

\begin{lemma}
\label{lem:GP}
Let Assumptions \ref{ass:convex} and \ref{ass:graph} hold.
If $\gamma_2 \ge (\overline l_F^2 / \underline l_F + \overline l_F)/\lambda_{\min}(\mathcal L)$, then $\mathbf{p}(t)$ approaches $\mathbf{p}^*$ with an exponential rate, where
$\mathbf{p}^* = 1_N \otimes p^*$, $p^*$ is the NE of $\mathbf{G}$,
and $\lambda_{\min}(\mathcal L)$ is the second minimal eigenvalue of $\mathcal L$.
\end{lemma}

\emph{Proof:}
Construct a Lyapunov function candidate $V_p$ as
$V_p(\mathbf{p}) = \frac 1 2 \Vert \mathbf{p}  - \mathbf{p}^*\Vert^2$.
With a similar procedure as the proof of Theorem $2$ in \cite{gadjov2018passivity}, there exists $\beta_0 > 0$ such that
\begin{equation}
\label{lya:grad-play:dev}
\dot V_p \le - \beta_0 \gamma_1 \Vert \mathbf{p}  - \mathbf{p}^*\Vert^2.
\end{equation}
Then $\Vert \mathbf{p}(t)  - \mathbf{p}^*\Vert^2 \le \exp{(-2\beta_0 \gamma_1)}\Vert \mathbf{p}(0)  - \mathbf{p}^*\Vert^2$, and the conclusion follows.
$\hfill\blacksquare$

\begin{remark}
Dynamics (\ref{gradient-play}) is inspired by (17) in  \cite{gadjov2018passivity}. Herein, $\gamma_1$ is used to adjust the convergence rate, and $\gamma_2$ relaxes the requirement on $\mathcal G$ since we do not impose assumptions on $\lambda_{\min}(\mathcal L)$. Note that the reference signal generator is introduced to compute an NE, and thus, other dynamics including the best-response and fictitious-play can be employed in place of the gradient-play scheme \cite{hu2022distributed, ye2023distributed}.
\end{remark}

\subsection{Internal Model}

In order to handle the external disturbances in (\ref{agent:dyn}), we design internal models, resorting to the ideas in \cite{huang2004general}. To begin with, we make the following assumptions.
\begin{assumption}
\label{ass:exosys}
The exosystem is neutrally stable, i.e., all eigenvalues of $S$ are semi-simple with zero real parts.
\end{assumption} 

\begin{assumption}
\label{ass:zero:state}
For each $i \in \mathcal I$, 
there exists a sufficiently smooth function $\mathbf{z}_i(s, v, w)$ with $\mathbf{z}_i(0, \mathbf{0}, w) = \mathbf{0}$
such that for all ${\rm col}\{v, w\} \in \mathbb R^{n_v + n_w}$ and $s \in \mathbb R$,
$\big({\partial \mathbf{z}_i(s, v, w)}/{\partial v}\big) Sv = f_{0i}(\mathbf{z}_i(s, v, w), s, v, w).$
\end{assumption}

Under Assumption \ref{ass:exosys}, given any compact set $\mathbb V_0 \subset \mathbb R^{n_v}$, there is a compact set $\mathbb V$ such that $v(t) \in \mathbb V, \forall t \ge 0$ if $v(0) \in \mathbb V_0$. Assumption \ref{ass:zero:state} is typical in solving  cooperative output regulation and optimal output consensus problems \cite{su2014cooperative, tang2020optimal, zhang2023distributed}.
Define
$\mathbf{z}_i^\star = \mathbf{z}_i(p_i^*, v, w),$
$\mathbf{x}_{1i}^\star = p_i$,
$\mathbf{x}_{2i}^\star =  - f_{1i}(\mathbf{z}_i^\star, p_i^*, v, w)$, $\dots$,
$\mathbf{x}_{(s+1)i}^\star = ({\partial \mathbf{x}_{si}^\star}/{\partial v}) Sv - f_{si}(\mathbf{z}_i^\star, p_i^*, \mathbf{x}_{2 i}^\star, \cdots, 
\mathbf{x}_{si}^\star, v, w)$, $s \in \{2, \cdots, r\}$ and 
$\mathbf{u}_i^\star = \mathbf{x}_{(r + 1)i}^\star$,
where $p^*$ is the NE of $\mathbf{G}$, and $p_i$ is given by (\ref{gradient-play}).

It is clear that 
$y_i(t)$ approaches $p_i(t)$ if there is a controller $u_i$ such that $\lim_{t \to \infty} \Vert x_{1i}(t) - \mathbf{x}_{1i}^\star \Vert = 0$.
In contrast to \cite{huang2004general, su2014cooperative, su2016cooperative}, the trajectory of $\mathbf{x}_{1i}^\star$ (or $p_i$) is generated by the reference signal generator (\ref{gradient-play}) instead of the exosystem (\ref{exosystem}), and as a result, it is more challenging to design internal models.
To facilitate the design, we adopt $p_i^*$ in the definition of  $\mathbf{x}_{(s+1)i}^\star$. Thus, $\mathbf{x}_{(s+1)i}^\star$ depends on $w$ as well as the unknown $p_i^*$. 

We make the following assumption for the existence of linear internal models.
\begin{assumption}
\label{ass:steady:state}
For each $i \in \mathcal I$, the functions $\mathbf{x}_{si}^\star, ~s\in \{2, \dots, r\}$ and $\mathbf{u}_i^\star$ are polynomials in $v$ with coefficients depending on $p_i^*$ and $w$.
\end{assumption}

In fact, Assumption \ref{ass:steady:state} is a well-suited condition, and has been widely used for the internal model design \cite{huang2004nonlinear}. It can be verified if $\mathbf{z}_i^\star$ and $f_{si}$ $s \in \{1, \cdots, r\}$ are all polynomials in their arguments $v, \mathbf{z}_i^\star, \mathbf{x}_{2 i}^\star, \cdots,  \mathbf{x}_{si}^\star$.
Under Assumption \ref{ass:steady:state}, there exist integers $n_s^i$ such that for any ${\rm col}\{v, w\} \in \mathbb V \times \mathbb W$ and $p^* \in \mathbb R^N$,
$\mathbf{x}_{(s+1)i}^\star$ satisfy
\begin{equation}
\begin{aligned}
\label{internal:ploy}
{{\rm d}^{n_s^i} \mathbf{x}_{(s+1)i}^\star}/{{\rm d}t^{n_s^i}} = \varsigma_{1i} \mathbf{x}_{(s+1)i}^\star + \varsigma_{2i}  {{\rm d} \mathbf{x}^\star_{(s+1)i}}/{{\rm d}t} + \cdots \\
+ \varsigma_{n_s^i i} {{\rm d}^{(n_s^i - 1)} \mathbf{x}^\star_{(s+1)i}} /{{\rm d}t^{(n_s^i - 1)}},
~s \in \{1, \dots, r\}
\end{aligned}
\end{equation}
where $\varsigma_{1i}, \dots, \varsigma_{n_s^i i}$ are scalars such that the roots of polynomials $P_s^i(\lambda) = \lambda^{n_s^i} - \varsigma_{1i} - \varsigma_{2i} \lambda - \cdots - \varsigma_{n_s^i i} \lambda^{n_s^i - 1}$ are distinct with zero real parts.
We should mention that the scalars $\varsigma_{1i}, \dots, \varsigma_{n_s^i i}$  are independent of $v$, $w$ and $p^*$.
Details for (\ref{internal:ploy}) can be found in \cite[Chapter 6]{huang2004nonlinear}, \cite{huang2004general}.
Define
\begin{equation}
\begin{aligned}
\label{def:Phi:Gamma}
\Phi_{si} = \left[\begin{array}{c|c}
0_{n_s^i -1} &I_{n_s^i -1} \\
\hline
\varsigma_{1i}  & \varsigma_{2i}, \dots, \varsigma_{n_s^i i}
\end{array}\right]\!, ~
\Gamma_{si}= \big[1, 0_{n_s^i -1}^{\top} \big].
\end{aligned}
\end{equation}
Let $M_{si} \in \mathbb R^{n_s^i \times n_s^i}$ be a Hurwitz matrix, and $N_{si} \in \mathbb R^{n_s^i}$ be a  vector such that the pair $(M_{si}, N_{si} )$ is controllable.
Since $(\Gamma_{si}, \Phi_{si})$ is observable, there is a nonsingular matrix $T_{si}$ satisfying 
$T_{si} \Phi_{si} - M_{si} T_{si} =  N_{si} \Gamma_{si}.$
Take
$\theta_{si} = T_{si} {\rm col}\{\mathbf{x}_{(s+1)i}^\star,\\
{{\rm d} \mathbf{x}_{(s+1)i}^\star}/{{\rm d}t}, \dots, 
{{\rm d}^{n_s^i - 1} \mathbf{x}_{(s+1)i}^\star}/{{\rm d}t^{n_s^i - 1}} \}$.
Then 
\begin{equation}
\begin{aligned}
\label{steady:state:generator}
\dot\theta_{si} = T_{si} \Phi_{si}T_{si}^{-1} \theta_{si},~
\mathbf{x}_{(s+1)i}^\star = \Psi_{si} \theta_{si}
\end{aligned}
\end{equation}
where $\Psi_{si} = \Gamma_{si} T_{si}^{-1}$.
Thus, system (\ref{steady:state:generator}) can be employed to generate the states $\mathbf{x}_{(s+1) i}^\star, s \in \{1, \dots, r\}$.
On the basis, we further design internal models for (\ref{agent:dyn}) as
\begin{equation}
\begin{aligned}
\label{internal:model}
\dot \eta_{si} &= M_{si} \eta_{si} + N_{si} x_{(s+1)i}, ~s \in \{1, \dots, r-1 \} \\
\dot \eta_{ri} &= M_{r_i i} \eta_{r i}+N_{ri} u_i,~i \in \mathcal I.
\end{aligned}
\end{equation}

\subsection{Augmented System}

By combining (\ref{agent:dyn}), (\ref{gradient-play}) with (\ref{internal:model}), we obtain an augmented dynamical system, for which we perform the following coordinate transformation
\begin{equation}
\begin{aligned}
\label{transformation}
&\bar z_i = z_i - \mathbf{z}_i^\star, ~~
\bar x_{1i} = x_{1i} - \mathbf{x}_{1i}^\star \\
&\bar x_{(s+1)i} = x_{(s+1)i} - \Psi_{si} \eta_{si} \\
&\tilde{\eta}_{si} = \eta_{si} - \theta_{si} - N_{si} \bar x_{si}, 
~ s \in \{1, \dots, r\}
\end{aligned}
\end{equation}
where $x_{(r+1)i} = u_i$. Let $\bar u_i = \bar x_{(r+1)i}$.
For convenience, define ${\bar x}_{[s]i} = {\rm col}\{{\bar x}_{1i}, \dots, {\bar x}_{si}\}$,
${\tilde \eta}_{[s]i} = {\rm col}\{{\tilde \eta}_{1i}, \dots, {\tilde \eta}_{si}\}$, ${\mathbf{x}}_{[s]i}^\star = {\rm col}\{{\mathbf{x}}_{1i}^\star, \dots, {\mathbf{x}}_{si}^\star\}$
and $\mu = {\rm col}\{v, w\}$. 
As a result, the augmented system reads as
\begin{equation}
\begin{aligned}
\label{aug:sys}
\dot p_i &= - \gamma_1\nabla_i J_i(\mathbf{p}_i) - \gamma_1 \gamma_2\sum\nolimits_{j \in \mathcal I_i} a_{ij}(p_i - p^j_i) \\
\dot p^i_{k} &= - \gamma_1\gamma_2\sum\nolimits_{j \in \mathcal I_i} a_{ij}(p^i_k - p^j_k)\\
\dot {\bar z}_i &= \bar f_{0i}(\bar z_i, \bar x_{1i}, p_i, p_i^*, \mu) + \hat f_{0i}(p_i, p_i^*, \mu) \\
\dot {\tilde \eta}_{1i} &= M_{1i} \tilde \eta_{1i} 
- N_{1i} \hat f_{1i}(p_i, p_i^*, \mu) + N_{1i}  \dot p_i \\
&~~~~~~~~~~~~~~ + \bar \kappa_{1i}(\bar z_i, \bar x_{1i}, p_i, p_i^*, \mu)  \\
\dot {\bar x}_{1i} &= \bar f_{1i}(\bar z_i, \bar x_{1i}, \tilde \eta_{1i}, p_i, p_i^*, \mu) \!+\! \hat f_{1i}(p_i, p_i^*, \mu) \!+\! \bar x_{2i} 
\!-\! \dot p_i \\
\dot {\tilde \eta}_{si} &= M_{si} \tilde \eta_{si} 
- N_{si} \hat f_{si}(p_i, p_i^*, \mu) \\
&~~~~~~~~~~~~~~ + \bar \kappa_{si}(\bar z_i, \bar x_{[s]i}, \tilde \eta_{[s-1]i}, p_i, p_i^*, \mu) \\
\dot {\bar x}_{si} &= \bar f_{si}(\bar z_i, \bar x_{[s]i}, \tilde \eta_{[s]i}, p_i, p_i^*, \mu) \!+\! \hat f_{si}(p_i, p_i^*, \mu) \!+\! \bar x_{(s+1)i} \\
&~~~~~~~~~~~~~~~~~~~~~~~~~~~~~~~~~~~~~~ 
s \in \{2, \dots, r\}
\end{aligned}
\end{equation}
where 
\begin{equation*}
\begin{aligned}
&\bar f_{0i}(\bar z_i, \bar x_{1i}, p_i, p_i^*, \mu) = f_{0i}(\bar z_i + \mathbf{z}_i^*, \bar x_{1i} + \mathbf{x}_{1i}^\star, \mu)\\
&~~~~~~~~~~~~~~~~ -f_{0i}(\mathbf{z}_i^*, p_i, \mu) \\
&\hat f_{0i}(p_i, p_i^*, \mu) = f_{0i}(\mathbf{z}_i^*, p_i, \mu) - f_{0i}(\mathbf{z}_i^*, p_i^*, \mu) \\
&\bar f_{1i}(\bar z_i, \bar x_{1i}, \tilde \eta_{1i}, p_i, p_i^*, \mu) = f_{1i}(\bar z_i + \mathbf{z}_i^*, \bar x_{1i} + \mathbf{x}_{1i}^\star, \mu)\\
&~~~~~~~~~~~~~~~~ -f_{1i}(\mathbf{z}_i^*, p_i, \mu)
+ \Psi_{1i}(\tilde \eta_{1i} + N_{1i} \bar x_{1i}) \\
&\hat f_{1i}(p_i, p_i^*, \mu) = f_{1i}(\mathbf{z}_i^*, p_i, \mu) - f_{1i}(\mathbf{z}_i^*, p_i^*, \mu) \\
&\bar \kappa_{1i}(\bar z_i, \bar x_{1i}, p_i, p_i^*, \mu) = M_{1i} N_{1i} \bar x_{1i} \\
&~~ -\! N_{1i} f_{1i}(\bar z_i + \mathbf{z}_i^*, \bar x_{1i} + \mathbf{x}_{1i}^\star, \mu) 
+ N_{1i} f_{1i}(\mathbf{z}_i^*, p_i, \mu) \\
\end{aligned}
\end{equation*}
\begin{equation*}
\begin{aligned}
&\bar f_{si}(\bar z_i, \bar x_{[s]i},  \tilde \eta_{[s]i}, p_i, p_i^*, \mu) = \Psi_{si}(\tilde \eta_{si} + N_{si} \bar x_{si}) \\
&~~~~~~~~~~~~~~~~ + \bar \omega_{si}(\bar z_i, \bar x_{[s]i}, \tilde \eta_{[s-1]i}, 
p_i, p_i^*, \mu) \\
&\hat f_{si}(p_i, p_i^*, \mu) \!=\! f_{si}(\mathbf{z}_i^*, \mathbf{x}_{[s]i}^\star, \mu) \\
&~~~~~~~~~~~~~~~~ -f_{si}(\mathbf{z}_i^*, p_i^*, \mathbf{x}_{2i}^\star, \dots, \mathbf{x}_{si}^\star, \mu),\\
&\bar \kappa_{si}(\bar z_i, \bar x_{[s]i}, \tilde \eta_{[s-1]i}, p_i, p_i^*, \mu) = M_{si} N_{si} \bar x_{si}  \\
&~~~~~~~~~~~~~~~~ - N_{si}\bar \omega_{si}\big(\bar z_i, \bar x_{[s]i}, \tilde \eta_{[s-1]i}, p_i, p_i^*, \mu \big)
\end{aligned}
\end{equation*}
with
\begin{equation*}
\begin{aligned}
&\bar \omega_{si}\big(\bar z_i, \bar x_{[s]i}, \tilde \eta_{[s-1]i}, p_i,  p_i^*, \mu \big) = f_{si}\big(\bar z_i + \mathbf{z}_i^\star, \bar x_{1i} + \mathbf{x}_{1i}^\star,  \\
&~~~~\bar x_{2i} + \Psi_{1i}(\tilde \eta_{1i} + \theta_{1i} + N_{1i} \bar x_{1i}), \dots,  \\
&~~~~\bar x_{si} + \Psi_{(s-1)i}(\tilde \eta_{(s-1)i} \!+\! \theta_{(s-1)i} \!+\! N_{(s-1)i} \bar x_{(s-1)i}), \mu \big)\\
&- f_{si}(\mathbf{z}_i^*, \mathbf{x}_{[s]i}^*, \mu) + f_{si}(\mathbf{z}_i^*, p_i^*, \mathbf{x}_{2i}^*, \dots, \mathbf{x}_{si}^*, \mu) + \Psi_{si} \theta_{si} \\
&- ({\partial \mathbf{x}_{si}^*}/{\partial v})S v\!-\! \Psi_{(s-1)i} \dot {\tilde \eta}_{(s-1)i} \!-\! \Psi_{(s-1)i} N_{(s-1)i}\dot {\bar x}_{(s-1)i}.
\end{aligned}
\end{equation*}

\begin{remark}
\label{rmk:aug}
Compared with the augmented systems for tackling output regulation problems in \cite{su2014cooperative, su2015cooperative, su2016cooperative}, (\ref{aug:sys}) is more complicated due to the presence of the reference signal generator (\ref{gradient-play}) and the extra terms $\hat f_{si}$.
Indeed, $\hat f_{si}$ originate from the definition of $\mathbf{z}_i^*$ and $\mathbf{x}_{si}^*$. Thus, it is more difficult to design controllers to stabilize (\ref{aug:sys}).
\end{remark}

We consider a class of distributed state-feedback control laws as
\begin{equation}
\label{cont:state:feedback}
\bar u_i = \varphi_i \big(\bar x_{1i}, \bar x_{2i}, \dots, \bar x_{ri}\big)
\end{equation}
where $\varphi_i$ is a sufficient smooth function vanishing at the origin. 
Let $\bar x_c = {\rm col}\{\mathbf{p}, \bar z_1, \bar x_{11} \dots, \bar x_{r1}, \tilde \eta_{11}, \dots, \tilde \eta_{r1}, \dots,  \bar z_N, \\
\bar x_{1N} \dots, \bar x_{rN}, \tilde \eta_{1N}, \dots, \tilde \eta_{rN}\} \in \mathbb R^{\bar n_c}$, where 
$\bar n_c = \sum_{i \in \mathcal I}\big(N + n_{z_i} + r + \sum_{s = 1}^r n_s^i \big)$.
Clearly, (\ref{aug:sys}) and (\ref{cont:state:feedback}) form a closed-loop system, one of whose equilibria is $\bar x_c = {\rm col}\{\mathbf{p}^*, 0_{\bar n_c - N^2}\}$, where $\mathbf{p}^* = 1_N \otimes p^*$. We further formulate a semi-global stabilization problem as follows.

\begin{problem}
\label{prb:stability}
Given any real number $\bar R > 0$, and any compact set $\mathbb V \times \mathbb W \subset \mathbb R^{n_v + n_w}$ containing the origin,  find a controller of the form  (\ref{cont:state:feedback}) such that for all $\mu \in \mathbb V \times \mathbb W$, the equilibrium point $\bar x_c = {\rm col}\{\mathbf{p}^*, 0_{\bar n_c - N^2}\}$ of the closed-loop system, composed of (\ref{aug:sys}) and (\ref{cont:state:feedback}), is asymptotically stable with its domain of attraction containing 
$\bar {\mathbb Q}^{\bar n_c}_{\bar R}$.
\end{problem}

The following lemma addresses the relation between the solvability of Problems \ref{prb} and \ref{prb:stability}.

\begin{lemma}
\label{lem:pro:conversion}
Let Assumptions \ref{ass:convex}-\ref{ass:steady:state} hold.
For any real number $\bar R > 0$ and any compact set $\mathbb V \times \mathbb W \in \mathbb R^{n_v + n_w}$ containing the origin, if Problem \ref{prb:stability} is solvable by a controller of the form (\ref{cont:state:feedback}), then Problem \ref{prb} can be solved by
\begin{equation}
\begin{aligned}
\label{cont:state:feedback1}
u_i &= \varphi_i \big(x_{1i} - p_i, x_{2i} - \Psi_{1i} \eta_{1i}, \dots, \\ &~~~~~~~~~x_{ri} - \Psi_{(r-1)i} \eta_{(r-1)i}\big) 
+ \Psi_{ri} \eta_{ri} \\
\dot \eta_{1i} &= M_{1i} \eta_{1i} + N_{1i} x_{2i} \\
& \vdots \\
\dot \eta_{ri} &= M_{r_i i} \eta_{r i} + N_{ri} u_i.
\end{aligned}
\end{equation}
\end{lemma}

\emph{Proof:}
For any $R > 0$ and $x_c(0) \in \bar{\mathbb Q}_R^{n_c}$, there exists $\bar R > 0$ such that $\bar x_c(0) \in \bar{\mathbb Q}_R^{n_c}$. 
By Assumption \ref{ass:exosys}, if $v(0) \in \mathbb V_0$, $v(t) \in \mathbb V, \forall t \ge 0$ for some compact set $\mathbb V$.
The solvability of Problem \ref{prb:stability} implies that for any
$\bar x_c(0) \in \bar {\mathbb Q}^{\bar n_c}_{\bar R}$,
the trajectory of $\bar x_c(t)$ is bounded for all $t \ge 0$, and moreover,  approaches ${\rm col}\{\mathbf{p}^*, 0_{\bar n_c - N^2}\}$.
Recalling (\ref{transformation}) yields the boundedness of $x_c(t)$ for all $t \ge 0$.
In addition, 
$\lim_{t \to \infty} \Vert y_i(t) - p_i^*\Vert 
\le \lim_{t \to \infty}   \Vert \bar x_{1i}(t) \Vert 
+ \lim_{t \to \infty}   \Vert p_i(t) - p_{i}^* \Vert$.
In light of Lemma \ref{lem:GP}, $\lim_{t \to \infty}   \Vert p_i(t) - p_{i}^* \Vert = 0$. It follows that $\lim_{t \to \infty} \Vert y_i(t) - p_i^*\Vert = 0$.
Thus, Problem \ref{prb} can be solved by (\ref{cont:state:feedback1}), and the proof is completed.
$\hfill\blacksquare$

\begin{remark}
Based on Lemma \ref{lem:pro:conversion}, our distributed NE seeking problem with dynamic agents is reformulated as a distributed semi-global stabilization problem of an augmented system.
We focus on the semi-global stabilization of (\ref{aug:sys}) since
we do not impose any Lipschitz conditions on the nonlinear functions $f_{si}$ in (\ref{agent:dyn}) as that of \cite{wang2015distributed}.
\end{remark}

\section{MAIN RESULTS}
\label{sec:results}

In this section, we solve Problem \ref{prb} by stabilizing the system (\ref{aug:sys}) semi-globally.
In order to do so, we make the following assumption on 
the zero dynamics of (\ref{aug:sys}).

\begin{assumption}
\label{ass:zero:dyn}
For each $i \in \mathcal I$, there exists a $C^2$ 
positive definite and proper function $V_{\bar z_i}: \mathbb R^{n_{z_i}} \rightarrow \mathbb R$ such that for all $\mu \in \mathbb V \times \mathbb W$ and $p_i \in \mathbb R$,
\begin{equation}
\label{ass:zero:lya}
\big(\partial V_{\bar z_i}(\bar z_i) / \partial \bar z_i \big) \cdot  \bar f_{0i}(\bar z_i, 0, p_i, p_i^*, \mu) \le  - \alpha_{i0}\Vert \bar z_i \Vert^2
\end{equation}
where $\alpha_{0i}$ is a known positive real number.
\end{assumption}

In fact, Assumption \ref{ass:zero:dyn} is quite standard. It implies that the zero dynamics of each agent is globally asymptotically stable as well as locally exponentially stable, and it is less stringent than the assumption of input-to-state stability in \cite{tang2020optimal}. A similar assumption can be found in   \cite{wang2015distributed, su2014cooperative, zhang2023distributed}.

Now we are ready to stabilize (\ref{aug:sys}) via a backstepping procedure. The main result is given as follows.

\begin{theorem}
\label{thm:semi-global:stability}
Let Assumptions \ref{ass:convex}-\ref{ass:zero:dyn} hold.
Given any real number $\bar R > 0$ and any compact set $\mathbb V \times \mathbb W \subset \mathbb R^{n_v + n_w}$ containing the origin, there exist $\gamma_1 > 0$ and $k_{si} > 0, s \in \{1, \dots, r\}, i \in \mathcal I$ depending on $\bar R$ such that for all $\mu \in \mathbb V \times \mathbb W$, the equilibrium point ${\rm col}\{\mathbf{p}^*, 0_{\bar n_c - N^2}\}$ of the augmented system (\ref{aug:sys}) is locally asymptotically stable with its region of attraction containing $\bar{\mathbb Q}_{\bar R}^{\bar n_c}$ by a distributed state-feedback controller as
\begin{equation}
\begin{aligned}
\label{cont:law}
\bar u_i = \!-\! \big[k_{ri} \bar x_{ri} \!+\! k_{ri} k_{(r-1)i}\bar x_{(r-1)i} \!+\! \cdots \!+\! \big(k_{ri} \cdots k_{1i}\big) \bar x_{1i}\big].
\end{aligned}
\end{equation}
\end{theorem}

The controller (\ref{cont:law}) takes the form of (\ref{cont:state:feedback}), and is determined by a recursive process. To start with, we first put (\ref{aug:sys}) into a block lower-triangular form.
Let 
\begin{equation}
\begin{aligned}
\label{coordinate:transf:x}
&\hat x_{1i} = \bar x_{1i} \\
&\hat x_{(s+1)i} = \bar x_{(s+1)i} + k_{si} \hat x_{si} \\
&\bar u_i = -k_{ri} {\hat x}_{ri}, ~ s \in \{1, \dots, r-1\}.
\end{aligned}
\end{equation}
Define
$Z_0 = {\rm col}\{\bar z_1, \dots, \bar z_N\}$,
$\chi_s = {\rm col}\{\tilde \eta_{s1}, \dots, \tilde \eta_{sN}\}$,
$\vartheta_s = {\rm col}\{\hat x_{s1}, \dots, \hat x_{sN}\}$, and
$K_s = {\rm diag}\{k_{s1}, \dots, k_{sN}\}$.
Then the system (\ref{aug:sys}) is cast into
\begin{equation}
\begin{aligned}
\label{aug:sys:sim}
\dot{\mathbf{p}} \!=&\! -\! \gamma_1 {\mathcal R}^{\top} \mathbf{F}(\mathbf{p})
\!-\! \gamma_1 \gamma_2 \mathbf{L} \mathbf{p} \\
\dot Z_0 \!=& H_0(Z_0, \vartheta_1, p, p^*, \mu) \!+\! \Upsilon_0(p, p^*, \mu) \\
\dot \chi_1 \!=& M_1 \chi_1 \!-\! N_1 \Upsilon_1(p, p^*\!, \mu) \!+\! N_1 \dot p \!+\! G_1(Z_0, \vartheta_1, p, p^*\!, \mu) \\
\dot {\vartheta}_1 \!=& H_1(Z_0, \vartheta_1,  \chi_1, p, p^*\!, \mu) \!+\! \Upsilon_1(p, p^*\!, \mu) \!-\! \dot p \!-\! K_1 \vartheta_1 \!+\! \vartheta_2 \\
\dot \chi_s  \!=& M_s \chi_s \!- \! N_s \Upsilon_s(p, p^*\!, \mu) \!+\!  G_s(Z_0, \vartheta_{[s]}, \chi_{[s-1]}, p, p^*\!, \mu) \\
\dot {\vartheta}_s \!=& H_s(Z_0,  \vartheta_{[s]}, \chi_{[s]}, p, p^*\!, \mu) \!+\! \Upsilon_s(p, p^*\!, \mu) \!-\! K_s \vartheta_s \!+\! \vartheta_{s + 1}
\end{aligned}
\end{equation}
where $H_0(\cdot) \!=\! {\rm col}\{\bar f_{01}(\cdot), \dots, \bar f_{0N}(\cdot)\}$,
$\Upsilon_0(\cdot) \!=\! {\rm col}\{\hat f_{01}(\cdot), \\
\dots, \hat f_{0N}(\cdot)\}$,
$M_1 = {\rm blkdiag}\{M_{11}, \dots, M_{1N}\} \in \mathbb R^{\tilde n_1 \times \tilde n_1}$,
$N_1 = {\rm blkdiag}\{N_{11}, \dots, N_{1N}\} \in \mathbb R^{\tilde n_1 \times N}$, $\tilde n_1 = \sum_{i \in \mathcal I} n_1^i$.
$H_1(\cdot) ={\rm col}\{\bar f_{11}(\cdot), \dots, \bar f_{1N}(\cdot)\}$,
$G_1(\cdot) = {\rm col}\{\bar \kappa_{11}(\cdot), \dots, \\
\bar \kappa_{1N}(\cdot)\}$, and
$\Upsilon_1(\cdot) \!=\! {\rm col}\{\hat f_{11}(\cdot), 
\dots, \hat f_{1N}(\cdot)\}$.
In addition, for all $s \in \{2, \dots, r\}$, 
$M_s = {\rm blkdiag}\{M_{s1}, \dots, M_{sN}\} \in \mathbb R^{\tilde n_s \times \tilde n_s}$,
$N_s = {\rm blkdiag}\{N_{s1}, \dots, N_{sN}\} \in \mathbb R^{\tilde n_s \times N}$,
$\tilde n_s = \sum_{i \in \mathcal I} n_s^i$,
$\hat x_{[s]i} = {\rm col}\{\hat x_{1i}, \dots, \hat x_{si}\}$,
$\chi_{[s]} = {\rm col}\{\chi_1, \dots, 
\chi_s\}$, 
$\vartheta_{[s]} = {\rm col}\{\vartheta_1, \dots, \vartheta_s\}$,
\begin{equation*}
\begin{aligned}
H_s(\cdot) \!&={\rm col}\{\bar f_{s1}(\bar z_1,
\hat x_{[s]1}, \tilde \eta_{[s]1}, p_1, p_1^*\!, \mu) \!+\! k_{(s-1)1} \dot{\hat x}_{(s-1)1}, \\
\dots,& ~\bar f_{sN}(\bar z_N,
\hat x_{[s]N}, \tilde \eta_{[s]N}, p_N, p_N^*\!, \mu) \!+\! k_{(s-1)N} \dot{\hat x}_{(s-1)N}\} \\
G_s(\cdot) 
\!&= {\rm col}\{\bar \kappa_{s1}(\bar z_1, \hat x_{[s]1}, \tilde \eta_{[s-1]1}, p_1, p_1^*\!, \mu), \dots, \\
&~~~~~~~~~~\bar \kappa_{sN}(\bar z_N, \hat x_{[s]N}, \tilde \eta_{[s-1]N}, p_N, p_N^*, \mu)\} \\
\Upsilon_s(\cdot) \!&= {\rm col}\{\hat f_{s1}(p_1, p_1^*, \mu), \dots, \hat f_{sN}(p_N, p_N^*, \mu)\}
\end{aligned}
\end{equation*}
with
\begin{equation*}
\begin{aligned}
&\bar f_{si}(\bar z_i, \hat x_{[s]i}, \tilde \eta_{[s]i}, p_i, p_i^*, \mu) = 
\bar f_{si}(\bar z_i, \hat x_{2i} - k_{1i} \hat x_{1i},  \\
&~~~~~~ \dots, \hat x_{si} - k_{si} \hat x_{(s-1)i}, \tilde \eta_{[s]i}, p_i, p_i^*, \mu)\\
&\bar \kappa_{si}(\bar z_i, \hat x_{[s]i}, \tilde \eta_{[s-1]i}, p_i, p_i^*, \mu)=\bar \kappa_{si}(\bar z_i, \hat x_{1i}, \hat x_{2i} - k_{1i} \hat x_{1i}, \\
&~~~~~~ \dots, \hat x_{si} - k_{si} \hat x_{(s-1)i}, \tilde \eta_{[s-1]i}, p_i, p_i^*, \mu),  i \in \mathcal I.
\end{aligned}
\end{equation*}
\emph{Proof of Theorem \ref{thm:semi-global:stability}:}
The proof is divided into the following three steps.

\emph{Step 1:} Analyze ${\rm col}\{\mathbf{p}, Z_0, \chi_1\}$-subsystem with $\vartheta_1 = \mathbf{0}$.

For the $\mathbf{p}$-subsystem, recalling (\ref{lya:grad-play:dev}) gives 
$\dot V_p \le - \beta_0 \gamma_1 \Vert \tilde{\mathbf{p}} \Vert^2$,
where $\tilde{\mathbf{p}} = \mathbf{p} - \mathbf{p}^*$.
For the $Z_0$-subsystem, define $V_{z}(Z_0) = \sum\nolimits_{i \in \mathcal I} V_{\bar z_i}(\bar z_i)$. By Assumption \ref{ass:zero:dyn}, we obtain
\begin{equation*}
\begin{aligned}
\dot V_z &= {\partial V_z}/{\partial Z_0} \big[H_0(Z_0, \mathbf{0}, p, p^*, \mu) + \Upsilon_0(p, p^*, \mu)\big]  \\
&\le  -\alpha_0\Vert Z_0\Vert^2 + \Vert {\partial V_z}/{\partial Z_0} \Vert \!\cdot\! \Vert \Upsilon_0(p, p^*, \mu) \Vert.
\end{aligned}
\end{equation*}
where $\alpha_0 = \min\{\alpha_{0i}\}_{i \in \mathcal I}$.
For the $\chi_1$-subsystem, there exists 
a positive define matrix
$P_1 \in \mathbb R^{\tilde n_1 \times \tilde n_1}$ such that $M_1^{\top} P_1 + P_1 M_1 \le - I_{\tilde n_1}$
since $M_1$ is a Hurwitz matrix.
Let $V_{\chi_1}(\chi_1) = \chi_1^{\top} P_1 \chi_1$. As a result,
\begin{equation*}
\begin{aligned}
\dot V_{\chi_1} &\le - \Vert \chi_1 \Vert^2 
+ 2 \Vert \chi_1 \Vert \!\cdot\! \Vert P_1 N_1 \Upsilon_1(p, p^*, \mu) \Vert \\
&+ 2 \Vert \chi_1 \Vert \!\cdot \!\Vert P_1 N_1 \dot p\Vert 
+  2 \Vert \chi_1 \Vert\! \cdot \! \Vert P_1 G_1(Z_0, 0, p, p^*, \mu)\Vert.
\end{aligned}
\end{equation*}
Note that $\bar x_c \in \bar{\mathbb Q}_{\bar R}^{\bar n_c}$ implies $\{\mathbf{p}, Z_0, \chi_1\} \in \bar{\mathbb Q}_{\bar R_1}^{N^2 + n_z + \tilde n_1}$,
where $\bar R_1 = \bar R$ and $n_z = \sum_{i \in \mathcal I} n_{z_i}$. 
By definition, there is $c_1 > 0$ such that
$\bar {\mathbb Q}_{\bar R_1}^{(N^2 + n_z + \tilde n_1)} \subset \bar{\Omega}_{c_1}(V_p) \times \bar{\Omega}_{c_1}(V_z) \times \bar{\Omega}_{c_1}(V_{\chi_1})$.

Construct a function as
\begin{equation}
\begin{aligned}
\label{pf:U1:def}
U_1(\mathbf{p}, Z_0, \chi_1) = V_p + \frac {c_1 V_z}{c_1 + 1 - V_z} + \zeta_1 V_{\chi_1}
\end{aligned}
\end{equation}
where $\zeta_1 > 0$. Clearly, $U_1$ is positive definite on $\mathbb R^{N^2} \times \Omega_{c_1 + 1}(V_z) \times \mathbb R^{\tilde n_1}$.
Define $\iota_1 = c_1^2 + (1 + \zeta_1)c_1$.
By \cite[Lemma 6]{su2014cooperative}, 
$\bar \Omega_{c_1}(V_p) \times  \bar \Omega_{c_1}(V_z) \times \bar \Omega_{c_1}(V_{\chi_1}) \subset \bar \Omega_{\iota_1}(U_1)$, and
$\bar \Omega_{\iota_1 + 1}(U_1) \subset \bar \Omega_{\iota_1+1} (V_p) \times \Omega_{c_1 + 1} (V_z)\times \bar\Omega_{\iota_1 + 1}(\zeta_1 V_z)$.
Besides, for any ${\rm col}\{\mathbf{p}, Z_0, \chi_1\} \in \bar \Omega_{\iota_1 + 1}(U_1)$, it holds that
$\bar L_{11} \le  {c_1 (c_1 + 1)}/{(c_1 + 1 - V_z)^2} \le 
\bar L_{12}$,
where $\bar L_{11} = {c_1}/{(c_1 + 1)}$ and
$\bar L_{12} = {(c_1 + \iota_1 + 1)^2}/(c_1^2 + c_1)$.
In light of \cite[Lemma 2]{su2014cooperative}, 
for all
$\mu \in \mathbb V \times \mathbb W$ and
${\rm col}\{\mathbf{p}, Z_0, \chi_1\} \in \bar \Omega_{\iota_1+1}(U_1)$,
there exist constants $\rho_{Z_0}$, $\rho_{\Upsilon_0}$, $\rho_p$, $\rho_{\Upsilon_1}$ and $\rho_{G_1}$ such that
$\Vert {\partial V_{z}}/{\partial Z_0}\Vert \le \rho_{Z_0} \Vert Z_0 \Vert$,
$\Vert \Upsilon_0(p, p^*, \mu) \Vert \le \rho_{\Upsilon_0} \Vert \tilde{\mathbf{p}}\Vert$,
$\Vert P_1 N_1 \dot p\Vert \le \rho_p \Vert \tilde{\mathbf{p}}\Vert$, 
$\Vert P_1 N_1 \Upsilon_1(p, p^*\!, \mu) \Vert \!\le\!   \rho_{\Upsilon_1} \Vert \tilde{\mathbf{p}}\Vert$, and
$\Vert  P_1 G_1(Z_0, 0, p,  p^*\!, \mu)\Vert
\le \rho_{G_1} \Vert Z_0 \Vert$.
Then
\begin{equation*}
\begin{aligned}
\dot U_1 =&~ \dot V_p + \frac {(c_1 + 1)c_1}{(c_1 + 1 - V_z)^2}  \dot V_z + \zeta_1 \dot V_{\chi} \\
\le&  -\beta_0 \gamma_1\Vert \tilde{\mathbf{p}} \Vert^2 
- \alpha_0 \bar L_{11} \Vert Z_0\Vert^2 
+ \rho_{Z_0} \rho_{\Upsilon_0}\bar L_{12} \Vert Z_0 \Vert \!\cdot\! \Vert \tilde{\mathbf{p}}\Vert \\
& - \zeta_1 \Vert \chi_1 \Vert^2 
+ \zeta_1\bar L_{13} \Vert \chi_1 \Vert \!\cdot\! \Vert \tilde{\mathbf{p}}\Vert 
+ \zeta_1\bar L_{14} \Vert \chi_1 \Vert \!\cdot\! \Vert Z_0 \Vert \\
\le& - \big(\beta_0 \gamma_1 - {\rho_{Z_0}^2 \rho_{\Upsilon_0}^2 \bar L_{12}^2}/(\alpha_0 \bar L_{11}) - \zeta_1\bar L_{13}^2) \Vert \tilde{\mathbf{p}}\Vert^2 \\
& - \big(3\alpha_0 \bar L_{11} / 4-  \zeta_1\bar L_{14}^2\big)
\Vert Z_0 \Vert^2 - \zeta_1/2 \cdot \Vert \chi_1 \Vert^2 
\end{aligned}
\end{equation*}
where $\bar L_{13} = 2  (\rho_p + \rho_{\Upsilon1})$, 
and $\bar L_{14} = 2 \rho_{G_1}$.
Take 
$\gamma_1 = {2 \rho_{Z_0}^2 \rho_{\Upsilon_0}^2 \bar L_{12}^2}/(\alpha_0\beta_0 \bar L_{11}) + 2\zeta_1\bar L_{13}^2 / \beta_0$,
$\zeta_1 = \alpha_0 \bar L_{11} / (4 \bar L_{14}^2)$
and
$\beta_1 = \min\{\beta_0 \gamma_1/2, \alpha_0 \bar L_{11} / 2 , \zeta_1/2\}$. For all
$\mu \in \mathbb V \times \mathbb W$ and
${\rm col}\{\mathbf{p}, Z_0, \chi_1\} \in \bar \Omega_{\iota_1+1}(U_1)$,
it holds that
\begin{equation}
\label{pf:U1:dev}
\dot U_1 \le - \beta_1 \big\Vert [\tilde{\mathbf{p}}; Z_0; \chi_1] \big\Vert^2.
\end{equation}

\emph{Step 2:}
Analyze the ${\rm col}\{\mathbf{p}, Z_0, \chi_1,  \vartheta_1\}$-subsystem with $\vartheta_2 = 0$.

Define
$V_{\vartheta_1}({\vartheta_1}) = \frac 1 2 \Vert \vartheta_1 \Vert^2$, and $k_1 = \min\{k_{1i}\}_{i \in \mathcal I}$. 
Then
\begin{equation*}
\begin{aligned}
\dot V_{\vartheta_1} 
\le& - k_1 \Vert \vartheta_1\Vert^2 +  \Vert \vartheta_1 \Vert \!\cdot\! \Vert H_1(Z_0, \vartheta_1,  \chi_1, p, p^*, \mu) \Vert \\
&+  \Vert \vartheta_1 \Vert \!\cdot\! \Vert \Upsilon_1(p, p^*, \mu)\Vert + \Vert \vartheta_1 \Vert \!\cdot\! \Vert \dot p\Vert.
\end{aligned}
\end{equation*}
Due to $\vartheta_1 \in \bar{\mathbb Q}^N_{\bar R_1}$, 
$V_{\vartheta_1} \le \hat c_1$ for some $\hat c_1 > 0$.
In the following, we consider the ${\rm col}\{\mathbf{p}, Z_0, \chi_1,  \vartheta_1\}$-subsystem with $\vartheta_2 = 0$.
Construct a function $W_1$ as
\begin{equation}
\begin{aligned}
\label{pf:W1:def}
W_1(\mathbf{p}, Z_0, \chi_1, \vartheta_1) = \frac {\iota_1 U_1}{\iota_1 + 1 - U_1} + \frac {\hat c_1 V_{\vartheta_1}} {\hat c_1 + 1 - V_{\vartheta_1}}.
\end{aligned}
\end{equation}
Then $W_1$ is positive definite on 
$\Omega_{\iota_1 + 1}(U_1) \times \Omega_{\hat c_1 + 1}(V_{ \vartheta_1})$. 
By \cite[Lemma 3]{su2016cooperative}, 
$\bar \Omega_{\iota_1}(U_1) \times \bar \Omega_{\hat c_1}(V_{ \vartheta_1}) \subset \bar \Omega_{\tau_1}(W_1)$, and
$\bar \Omega_{\tau_1 + 1}(W_1) \subset \Omega_{\iota_1 + 1}(U_1) \times \Omega_{\hat c_1 + 1}(V_{\vartheta_1})$,
where $\tau_1 = \iota_1^2 + \hat c_1^2$.
Moreover, for all ${\rm col}\{\mathbf{p}, Z_0, \chi_1, \vartheta_1\} \in \bar \Omega_{\tau_1 + 1}(W_1)$,
$\hat L_{11} \le {\iota_1(\iota_1 + 1)}/{(\iota_1 + 1 - U_1)^2} \le  \hat L_{12},$
and 
$\hat L_{13} \le {\hat c_1(\hat c_1 + 1)}/{(\hat c_1 + 1 - V_{\vartheta_1})^2} \le \hat L_{14}$
where $\hat L_{11} = {\iota_1}/{(\iota_1 + 1)}$,
$\hat L_{12} ={(\iota_1 + \tau_1 + 1)^2}/{(\iota_1^2 + \iota_1)}$,
$\hat L_{13} = {\hat c_1}/{(\hat c_1 + 1)}$, and
$\hat L_{14} = {(\hat c_1 + \tau_1 + 1)^2}/{(\hat c_1^2 + \hat c_1)}$.
Notice that we currently do not impose $\vartheta_1 = \mathbf{0}$.
Combining (\ref{aug:sys:sim}), (\ref{pf:U1:def}) with (\ref{pf:U1:dev}), we obtain
\begin{equation*}
\begin{aligned}
& \dot U_1 \le - \beta_1 \big\Vert [\tilde{\mathbf{p}}; Z_0; \chi_1] \big\Vert^2 \\
&\!+\! \Vert \partial U_1 / \partial Z_0 \Vert \!\cdot\! \Vert H_0(Z_0, \vartheta_1, p, p^*, \mu) \!-\! H_0(Z_0, \mathbf{0}, p, p^*, \mu)\Vert \\
&\!+\!  \Vert \partial U_1 / \partial \chi_1 \Vert \!\cdot\!
\Vert G_1(Z_0, \vartheta_1, p, p^*, \mu) - G_1(Z_0, \mathbf{0}, p, p^*, \mu) \Vert.
\end{aligned}
\end{equation*}
In light of \cite[Lemma 2]{su2014cooperative},
there exist positive real numbers $\sigma_{U_1}$, $\sigma_{H_0}$, $\sigma_{G_1}$, $\sigma_{H_1}$, $\sigma_{\Upsilon_1}$ and $\tilde \rho_p$ such that 
for all $\mu \in \mathbb V \times \mathbb W$ and ${\rm col}\{\mathbf{p}, Z_0, \chi_1, \vartheta_1\} \in \bar \Omega_{\tau_1 + 1}(W_1)$,
\begin{equation}
\begin{aligned}
\label{pf:ieqs}
&\Vert {\partial U_1}/{\partial Z_0}\Vert  \le \sigma_{U_1} \Vert Z_0 \Vert, ~
\Vert {\partial U_1}/{\partial \chi_1} \Vert  \le \sigma_{U_1} \Vert \chi_1 \Vert \\
&\Vert H_0(Z_0, \vartheta_1, p, p^*, \mu) - H_0(Z_0, \mathbf{0}, p, p^*, \mu) \Vert \le \sigma_{H_0} \Vert \vartheta_1 \Vert \\
&\Vert G_1(Z_0, \vartheta_1, p, p^*, \mu) - G_1(Z_0, \mathbf{0}, p, p^*, \mu) \Vert \le \sigma_{G_1} \Vert \vartheta_1 \Vert \\
&\Vert H_1(Z_0, \vartheta_1,  \chi_1, p, \mu) \Vert \le
\sigma_{H_1} \big(\Vert Z_0 \Vert + \Vert \vartheta_1 \Vert + \Vert \chi_1 \Vert\big)\\
&\Vert \Upsilon_1(p, p^*, \mu)\Vert  \le \sigma_{\Upsilon_1} \Vert \tilde{\mathbf{p}}\Vert,~
\Vert \dot p\Vert \le \hat \rho_p \Vert \tilde{\mathbf{p}}\Vert.
\end{aligned}
\end{equation}
It follows that
\begin{equation*}
\begin{aligned}
\dot W_1 =&~  \frac {\iota_1(\iota_1 + 1)}{(\iota_1 + 1 - U_1)^2} \dot U_1 
+ \frac {\hat c_1(\hat c_1 + 1)}{(\hat c_1 + 1 - V_{\vartheta_1})^2} 
\dot V_{\vartheta_1} \\
\le & - \beta_1 \hat L_{11} \big\Vert [\tilde{\mathbf{p}}; Z_0; \chi_1] \big\Vert^2 
+ \sigma_{U_1}\sigma_{H_0} \hat L_{12} \Vert Z_0\Vert \!\cdot \! \Vert \vartheta_1 \Vert \\
&+ \sigma_{U_1}\sigma_{G_1}\hat L_{12} \Vert \chi_1\Vert \!\cdot \! \Vert \vartheta_1 \Vert 
- k_1 \hat L_{13} \Vert \vartheta_1\Vert^2 \\
&+ \sigma_{H_1} \hat L_{14} \Vert \vartheta_1 \Vert  \big(\Vert Z_0 \Vert + \Vert \vartheta_1 \Vert + \Vert \chi_1 \Vert\big) \\
&+ \rho_{\Upsilon_1} \hat L_{14} \Vert \vartheta_1 \Vert \!\cdot\! \Vert \tilde{\mathbf{p}}\Vert
+  \hat \rho_p \hat L_{14} \Vert \vartheta_1 \Vert \!\cdot\! \Vert \tilde{\mathbf{p}}\Vert \\
\le&- \beta_1 \hat L_{11}/2 \cdot \big\Vert [\tilde{\mathbf{p}}; Z_0; \chi_1] \big\Vert^2 \\
&- \big(k_1 \hat L_{13} -\delta_1 - \sigma_{H_1} \hat L_{14} \big)\Vert \vartheta_1\Vert^2.
\end{aligned}
\end{equation*}
where
$\delta_1 = \big(\sigma_{U_1}^2 \hat L_{12}^2 ( \sigma_{H_0}^2 +  \sigma_{G_1}^2) +  \hat L_{14}^2 (\hat \rho_p^2 + \rho_{\Upsilon_1}^2) + 2\sigma_{H_1}^2 \hat L_{14}^2 \big)/(\beta_1 \hat L_{11})$. Take $k_1 = 2(\delta_1 +  \sigma_{H_1} \hat L_{14}) / \hat L_{13}$, and 
$\alpha_1 = \min\{\beta_1 \hat L_{11}/2, k_1 \hat L_{13}/2\}$.
Then, for all $\mu \in \mathbb V \times \mathbb W$ and ${\rm col}\{\mathbf{p}, Z_0, \chi_1, \vartheta_1\} \in \bar \Omega_{\tau_1 + 1}(W_1)$,
\begin{equation}
\label{pf:W1:dev}
\dot W_1 \le - \alpha_1 \Vert [\tilde{\mathbf{p}}; Z_0; \chi_1;  \vartheta_1] \Vert^2.
\end{equation}

\emph{Step 3:}
Prove Theorem $1$ by induction.

Define 
$X_1 = {\rm col}\{Z_0, \chi_1, \vartheta_1\}$, 
$\hat n_1 = N^2 + n_z + \tilde n_1 + N$,
$X_s = {\rm col}\{X_{s-1}, \chi_s, \vartheta_s\}$, and
$\hat n_s = \hat n_{s-1} + \tilde n_s + N$
for $s \in \{2, \dots, r\}$.
By (\ref{coordinate:transf:x}), $\bar x_c \in \bar{\mathbb Q}_{\bar R}^{\bar n_c}$ implies $\{\mathbf{p}, X_s\} \in \mathbb D_s \triangleq \bar{\mathbb Q}_{\bar R_1}^{\hat n_1} \times \bar{\mathbb Q}_{\bar R_2}^{\hat n_2 - \hat n_1} \times \cdots \times \bar{\mathbb Q}_{\bar R_s}^{\hat n_s - \hat n_{s-1}}$ for some positive real numbers $\bar R_1, \dots, \bar R_s$.

Consider the $\{\mathbf{p}, X_{s-1}\}$-subsystem with $ \vartheta_{s} = 0$. 
Based on (\ref{pf:W1:dev}), we suppose there is a continuously differentiable and positive definite function  $W_{s-1}(\mathbf{p}, X_{s-1})$ on  $\bar \Omega_{\tau_{s-1}+1}(W_{s-1})$
such that $\mathbb D_{s-1} \subset \bar \Omega_{\tau_{s-1}}(W_{s-1})$ for some $\tau_{s-1} > 0$.
Besides, for all $\mu \in \mathbb V \times \mathbb W$ and 
${\rm col}\{\mathbf{p}, X_{s-1}\} \in \bar \Omega_{\tau_{s-1} + 1} (W_{s-1})$, there is $\alpha_{s-1} > 0$ such that
$$\dot W_{s-1} \le - \alpha_{s-1} \big\Vert [\tilde{\mathbf{p}}; X_{s-1}] \big\Vert^2.$$

Recalling (\ref{aug:sys:sim}) gives
\begin{equation*}
\begin{aligned}
\dot \chi_s & = M_s \chi_s \!-\!  N_s \Upsilon_s(p, p^*, \mu)
\!+\! G_s(X_{s-1}, \vartheta_s, p, p^*, \mu),\\
\dot {\vartheta}_s & =  H_s(X_{s-\!1},\! \vartheta_s,\!
\chi_s, p, p^*\!,  \mu) \!+\! \Upsilon_s(p, p^*\!, \mu) \!-\! K_s  \vartheta_s \!+\! \vartheta_{s+1}.
\end{aligned}
\end{equation*}

It is clear that $M_s^{\top} P_s + P_s M_s \le - I_{\tilde n_s}$ for some 
positive definite matrix $P_s$. Define $V_{\chi_s}(\chi_s) = \chi_s^{\top} P_s \chi_s$. 
Due to $\chi_s \in \bar{\mathbb Q}_{\bar R_s}^{\tilde n_s}$, there is 
$c_s > 0$ such that $V_{\chi_s} \le c_s$.
Besides,
\begin{equation*}
\begin{aligned}
\dot V_{\chi_s} \le& - \Vert \chi_s \Vert^2 
+ 2\Vert  \chi_s \Vert \!\cdot\! \Vert P_s N_s \Upsilon_s(p, p^*, \mu) \Vert\\ 
&+  2\Vert  \chi_s \Vert \!\cdot\!  \Vert P_s G_s(X_{s-1},  \vartheta_s, p, p^*, \mu) \Vert.
\end{aligned}
\end{equation*}

Consider the $\{\mathbf{p}, X_{s-1}, \chi_{s}\}$-subsystem with 
$\vartheta_{s} = 0$. Let
\begin{equation*}
\begin{aligned}
U_s(\mathbf{p}, X_{s-1}, \chi_{s}) = \frac{\tau_{s-1} W_{s-1}}{\tau_{s-1} + 1 - W_{s-1}} + \zeta_{s} V_{\chi_s}.
\end{aligned}
\end{equation*}
Then $U_s$ is positive definite on $\Omega_{\tau_{s-1} + 1}(W_{s-1}) \times \mathbb R^{\tilde n_s}$.
$\bar \Omega_{\tau_{s-1}}(W_{s-1}) \times \bar \Omega_{c_s}(V_{\chi_s}) \subset \bar \Omega_{\iota_s}(U_s)$, and $\bar \Omega_{\iota_s + 1}(U_s) \subset \times \Omega_{\tau_{s-1} + 1}(W_{s-1}) \times \bar \Omega_{\iota_s + 1}(\zeta_{s} V_{\chi_s})$, 
where $\iota_s = \tau_{s-1}^2 + \zeta_s c_s$.
For all $\mu \in \mathbb V \times \mathbb W$ and ${\rm col}\{\mathbf{p}, X_{s-1}, \chi_s\} \in \bar \Omega_{\iota_s + 1}(U_s)$,
$$\bar L_{s1} \le  {\tau_{s-1} (\tau_{s-1} + 1)}/{(\tau_{s-1} + 1 - W_{s-1})^2} \le \bar L_{s2}
$$
where $\bar L_{s1} = {\tau_{s-1}}/{(\tau_{s-1} + 1)}$ and
$\bar L_{s2} = {(\tau_{s-1} + \iota_s + 1)^2}\\
/(\tau_{s-1}^2 + \tau_{s-1})$.
In addition,
$\Vert P_s N_s \Upsilon_s(p, p^*, \mu) \Vert \le \rho_{\Upsilon_s} \Vert \tilde{\mathbf{p}}\Vert$ and
$\Vert P_s G_s(X_{s-1}, \mathbf{0}, p, p^*, \mu) \Vert \le 
\rho_{G_s} \Vert  X_{s-1}\Vert$ for some 
$\rho_{\Upsilon_s}, \rho_{G_s} > 0$.
Hence,
\begin{equation*}
\begin{aligned}
\dot U_s =& ~ \frac{\tau_{s-1}(\tau_{s-1} + 1)} {(\tau_{s-1} + 1 - W_{s-1})^2} \dot W_{s-1} + \zeta_s \dot V_{\chi_s} \\
\le& - \alpha_{s-1} \bar L_{s1} \big\Vert [\tilde{\mathbf{p}}; X_{s-1}] \big\Vert^2 - \zeta_s \Vert \chi_s \Vert^2 \\
&+ 2 \zeta_s  \rho_{\Upsilon_s} \Vert  \chi_s \Vert \!\cdot\! \Vert \tilde{\mathbf{p}}\Vert +  2 \zeta_s  \rho_{G_s} \Vert \chi_s \Vert \!\cdot\!  \Vert X_{s-1} \Vert \\
\le & - \zeta_s /2 \cdot \Vert \chi_s \Vert^2 - \big(\alpha_{s-1} \bar L_{s1} - \hat \delta_s \zeta_s \big) \big\Vert [\tilde{\mathbf{p}}; X_{s-1}] \big\Vert^2
\end{aligned}
\end{equation*}
where $\hat \delta_s = 4 \zeta_s (\rho_{\Upsilon_s}^2 + \rho_{G_s}^2)$.
Let $\zeta_s = \alpha_{s-1} \bar L_{s1} / (2\hat \delta_s)$,
and $\beta_s = \min\{\zeta_s /2, \alpha_{s-1} \bar L_{s1}/2\}$. 
It follows that
\begin{equation*}
\begin{aligned}
\dot U_s \le -\beta_s \big\Vert [\tilde{\mathbf{p}}; X_{s-1}; \chi_s]\big\Vert^2.
\end{aligned}
\end{equation*}

Finally, we discuss the $\{\mathbf{p}, X_s\}$-subsystem with $\vartheta_{s+1} = \mathbf{0}$.
Define 
$V_{\vartheta_s} = \frac 12 \Vert \vartheta_s\Vert^2$ and $k_s = \min\{k_{si}\}_{i \in \mathcal I}$.
Then
\begin{equation*}
\begin{aligned}
\dot V_{\vartheta_s} 
\le& - k_s \Vert \vartheta_s\Vert^2 
+ \Vert \vartheta_s \Vert \!\cdot\! \Vert \Upsilon_s(p, p^*, \mu)\Vert\\
&+ \Vert \vartheta_s \Vert \!\cdot\! \Vert H_s(X_{s-1},  \vartheta_s, \chi_s, p, p^*, \mu) \Vert.
\end{aligned}
\end{equation*}

Since $\vartheta_s \in \bar{\mathbb Q}_{\bar R_s}^N$, 
$V_{\vartheta_s} \le \hat c_s$ for some $\hat c_s > 0$.
Let
\begin{equation*}
\begin{aligned}
W_s(\mathbf{p}, X_{s}) = \frac {\iota_s U_s}{\iota_s + 1 -U_s}  + \frac {\hat c_s V_{\vartheta_s}}{\hat c_s + 1 - V_{\vartheta_s}}.
\end{aligned}	
\end{equation*}
Then $W_s$ is positive definite on $\Omega_{\iota_s + 1}(U_s) \times \Omega_{\hat c_s + 1}(V_{\vartheta_s})$,
$\bar \Omega_{\iota_s}(U_s) \times \bar \Omega_{\hat c_s} (V_{\vartheta_s}) \subset \bar \Omega_{\tau_s}(W_s)$, and
$\bar \Omega_{\tau_s + 1}(W_s) \subset \Omega_{\iota_{s} + 1} (U_s) \times \Omega_{\hat c_{s} + 1} (V_{\vartheta_s})$,
where $\tau_s = \iota_s^2 + c_s^2$. 
For all $\mu \in \mathbb V \times \mathbb W$ and ${\rm col}\{\mathbf{p}, X_s\} \in \bar \Omega_{\tau_s + 1}(W_s)$,
$\hat L_{s1} \le {\iota_s(\iota_s + 1)}/{(\iota_s + 1 - U_s)^2} \le  \hat L_{s2}$,
and
$\hat L_{s3} \le {\hat c_s(\hat c_s \!+ 1)}/(\hat c_s \\
+ 1 - V_{\vartheta_s})^2 \le \hat L_{s4}$,
where $\hat L_{11} = {\iota_s}/{(\iota_s + 1)}$,
$\hat L_{s2} ={(\iota_s + \tau_s + 1)^2}/{(\iota_s^2 + \iota_s)}$,
$\hat L_{s3} = {\hat c_s}/{(\hat c_s + 1)}$, and
$\hat L_{s4} = {(\hat c_s + \tau_s + 1)^2}/{(\hat c_s^2 + \hat c_s)}$.
Notice that
\begin{equation*}
\begin{aligned}
&\dot U_s \le -\beta_s \big\Vert [\tilde{\mathbf{p}}; X_{s-1}; \chi_s]\big\Vert^2 + \Vert \partial U_s / \partial \vartheta_{s-1} \Vert \!\cdot\! \Vert \vartheta_s \Vert + \\
& \Vert \partial U_s / \partial \chi_s \Vert \!\cdot\! 
\Vert G_s(X_{s\!-\!1}, \vartheta_s, p, p^*\!, \mu) \!-\! G_s(X_{s\!-\!1}, \mathbf{0}, p, p^*, \mu) \Vert.
\end{aligned}	
\end{equation*}

Similar to (\ref{pf:ieqs}), there exist $\sigma_{U_s}$,
$\sigma_{H_s}$, $\sigma_{\Upsilon_s}$
and $\sigma_{G_s}$
such that for all 
$\mu \in \mathbb V \times \mathbb W$ and
${\rm col}\{\mathbf{p}, X_s\} \in \bar \Omega_{\tau_s + 1}(U_s)$, 
$\Vert \partial U_s /  \vartheta_{s-1}  \Vert \le \sigma_{U_s} \Vert  \vartheta_{s-1} \Vert$,
$\Vert \partial U_s /  \chi_s  \Vert \le \sigma_{U_s} \Vert \chi_s \Vert$,
$\Vert H_s(X_{s-1}, \vartheta_s, \chi_s, p, p^*, \mu) \Vert \le
\sigma_{H_s}(\Vert X_{s-1}\Vert + \Vert \vartheta_s\Vert + \Vert \chi_s\Vert)$,
$\Vert \Upsilon_s(p, p^*, \mu)\Vert \le \sigma_{\Upsilon_s} \Vert \tilde{\mathbf{p}}\Vert$ and
$\Vert G_s(X_{s-1}, \vartheta_s, p, p^*, 
\mu) - G_s(X_{s-1}, \mathbf{0}, p, p^*, \mu) \Vert \le \sigma_{G_s} \Vert  \vartheta_s \Vert$. 
It follows that
\begin{equation*}
\begin{aligned}
&\dot W_s \le -\beta_s \hat L_{s1}\big\Vert [\tilde{\mathbf{p}}; X_{s-1}; \chi_s]\big\Vert^2 - k_s \hat L_{s3} \Vert \vartheta_s\Vert^2 \\
&~+\! \sigma_{U_s} \hat L_{s2} \Vert \vartheta_s \Vert \!\cdot\!
\big( \Vert \vartheta_{s-1} \Vert \!+\!  \Vert  \chi_s \Vert\big)
\!+\! \sigma_{\Upsilon_s} \hat L_{s4} \Vert \vartheta_s \Vert \!\cdot\! 
\Vert  \tilde{\mathbf{p}}\Vert \\
&~+\! \sigma_{H_s} \hat L_{s4} \Vert \vartheta_s \Vert \! \cdot \!
\big(\Vert X_{s-1}\Vert \!+\! \Vert \vartheta_s\Vert + \Vert \chi_s\Vert \big) \\
& \le - \beta_s \hat L_{s1}/2 \cdot \big\Vert [\tilde{\mathbf{p}}; X_{s-1}; \chi_s]\big\Vert^2 -  (k_s \hat L_{s3} - \delta_s) 
\Vert \vartheta_s\Vert^2
\end{aligned}		
\end{equation*}	
where $\delta_s = \big(2\sigma_{U_s}^2 \hat L_{s2}^2 + (\sigma_{\Upsilon_s}^2 + 2\sigma_{H_s}^2)\hat L_{s4}^2 \big) / (\beta_s \hat L_{s1}) + \sigma_{H_s} \hat L_{s4}$.
Take $k_s = 2 \delta_s / \hat L_{s3}$, and 
$\alpha_s = \min \{\beta_s \hat L_{s1}/2, 
k_s \hat L_{s3}/2\}$.
Then
\begin{equation}
\begin{aligned}
\label{pf:Ws:dev}
\dot W_s \le - \alpha_s \big\Vert [\tilde{\mathbf{p}}, X_s] \big\Vert^2.
\end{aligned}
\end{equation}
When $s = r$ and $\vartheta_{r + 1} = \mathbf{0}$, it follows from (\ref{pf:Ws:dev}) that for all $\mu \in \mathbb V \times \mathbb W$ and ${\rm col}\{\mathbf{p}(0), X_r(0)\} \in \mathbb D_ r$, 
$\dot W_r \le - \alpha_r \big\Vert [\tilde{\mathbf{p}}, X_r] \big\Vert^2$ for some $\alpha_r > 0$.
Thus, the trajectory of ${\rm col}\{\mathbf{p}(t), X_r(t)\}$ is bounded, and moreover, converges to ${\rm col}\{\mathbf{p}^*, 0_{\hat n_r -  N^2}\}$. This completes the proof.
$\hfill\blacksquare$

\begin{remark}
In the above analysis, a backstepping procedure is employed since
dynamics (\ref{agent:dyn}) is in a lower-triangular form, which covers the systems in \cite{wang2015distributed, zhang2019distributed, tang2020optimal,zhang2023distributed}.
The proof of Theorem \ref{thm:semi-global:stability} is inspired by those of \cite{su2014cooperative, su2015cooperative, su2016cooperative}, but it is more challenging due to complexity of (\ref{aug:sys}) as mentioned in Remark \ref{rmk:aug}.
We overcome the obstacles by leveraging the exponential convergence of (\ref{gradient-play}), and carefully constructing the Lyapunov function candidates $U_s$ and $W_s$.
\end{remark}

By combining Lemma \ref{lem:pro:conversion} with Theorem \ref{thm:semi-global:stability}, we establish the following result.

\begin{theorem}
Let Assumptions \ref{ass:convex}-\ref{ass:zero:dyn} hold.
Given any $R > 0$ and any compact set $\mathbb V_0 \times \mathbb W \subset \mathbb R^{n_v + n_w}$, there exist $\gamma_1 > 0$ and $k_{si} > 0, s \in \{1, \dots, r\}, i \in \mathcal I$  depending on $ R$
such that Problem \ref{prb} is solvable by a distributed dynamic state-feedback controller as
\begin{equation*}
\begin{aligned}
\label{cotroller}
u_i =& - k_{ri}\big(x_{ri} - \Psi_{(r-1)i} \eta_{(r-1)i}\big) \\
&- k_{ri}k_{(r-1)i}\big(x_{(r-1)i} - \Psi_{(r-2)i} \eta_{(r-2)i}\big)
- \dots \\
&- k_{ri}k_{(r-1)i} \dots k_{2i}k_{1i}\big(x_{1i} -p_i\big) + \Psi_{ri} \eta_{ri} 
\end{aligned}
\end{equation*}

\begin{equation}
\begin{aligned}
\dot \eta_{1i} =&~ M_{1i} \eta_{1i} + N_{1i} x_{2i} \\
& \vdots \\
\dot \eta_{ri} =&~ M_{r_i i} \eta_{r i} + N_{ri} u_i.
\end{aligned}
\end{equation}
\end{theorem}

\begin{remark}
We summarize the procedure to solve Problem \ref{prb} as follows. First, construct a reference signal generator (\ref{gradient-play}) for distributed NE seeking, where the gains $\gamma_1$ and $\gamma_2$ can refer to Lemma \ref{lem:GP} and Theorem \ref{thm:semi-global:stability}. Second, find $\Phi_{si}$ and $\Gamma_{si}$ in (\ref{def:Phi:Gamma}), select $M_{si}$ and $N_{si}$ manually, and design internal models (\ref{internal:model}).
Third, derive the augmented system (\ref{aug:sys}).
Fourth, determine the controller (\ref{cotroller}), in which the gains $k_{si}$ are obtained by the recursive design given in the proof of Theorem \ref{thm:semi-global:stability}.
\end{remark}

\section{EXAMPLE}
\label{sec:example}

Consider a multi-agent system with four agents given by 
\begin{equation*}
\begin{aligned}
\dot z_i &= g_{1i} z_i + x_{1i} + g_{2i} v_1 \\
\dot x_{1i} &= g_{3i} z_i x_{1i} + g_{4i}v_2 + x_{2i} \\
\dot x_{2i} &= g_{5i} z_i^2 x_{1i} +  g_{6i}x_{1i} x_{2i} + u_i, 
i \in \{1, \dots, 4\} \\
\end{aligned}
\end{equation*}
where $g_i = {\rm col}\{g_{1i}, \dots, g_{6i}\} \in \mathbb R^6$ is an uncertain vector satisfying $g_{1i} < 0$.
The exosystem (\ref{exosystem}) is given by
$\dot v_1 = v_2$ and $\dot v_2 = - v_1$, where $v = [v_1, v_2]^{\top} \in \mathbb R^2$.
Besides, let agent $i$ be endowed with a local cost function as
$J_i(y_i, y_{-i}) = (y_i - h_{1i})^2 + y_i (h_{2i} \sum_{j \in \mathcal I} y_j + h_{3i})$, and all agents communicate over a ring graph, where $h_i = {\rm col}\{h_{1i}, h_{2i}, h_{3i}\} \in \mathbb R^3$ is known.
It is clear that Assumptions \ref{ass:convex} - \ref{ass:exosys} hold. 

By setting 
$\mathbf{z}_i(s, v, w) = -g_{1i}g_{2i} v_1 / (g_{1i}^2 + 1) - g_{2i} v_2 / (g_{1i}^2 + 1) - x_{1i} / g_{1i}$, 
Assumption \ref{ass:zero:state} is satisfied.
It follows that
$\mathbf{z}_i^\star = \mathbf{z}_i(p_i^*, v, w)$,
$\mathbf{x}_{1i}^\star = p_i$,
$\mathbf{x}_{2i}^\star = - g_{3i} p_i^* \mathbf{z}_i^\star - g_{4i}v_2$, 
$\mathbf{u}_i^\star = (\partial \mathbf{x}_{2i}^\star / \partial v) Sv - g_{5i} p_i^* \mathbf{z}_i^{\star, 2} -  g_{6i}p_i^* \mathbf{x}_{2i}^\star$, and $\bar f_{0i}(\bar z_i, 0, p_i, p_i^*, \mu) = g_{1i} \bar z_i$.
It is straightforward to verify that Assumptions \ref{ass:steady:state} and \ref{ass:zero:dyn} hold.
Notice that
${{\rm d}^3 \mathbf{x}_{2i}^\star}/{{\rm d} t^3} = - {{\rm d} \mathbf{x}_{2i}^\star}/{{\rm d} t}$, and
${{\rm d^5} \mathbf{u}_{i}^*}/{{\rm d} t^5} = -4 {{\rm d} \mathbf{u}_{i}^*}/{{\rm d} t} - 5 {{\rm d}^3 \mathbf{u}_{i}^*}/{{\rm d} t^3}$.
Then we obtain $\Phi_{si}, s \in \{1, 2\}$ via (\ref{def:Phi:Gamma}).
Given
\begin{equation*}
\begin{aligned}
M_{1i} \!=\! \left[\begin{array}{c|c}\!\!
\! 0_2\!\! & \!\! I_{2} \\
\hline
\!\!\!-3\!\!  &\!\! -7,\! -5\!\!\!
\end{array}\right]\!\!, 
M_{2i} \!=\! \left[\begin{array}{c|c}
\!\!\!0_4 \! & \!\! I_{4} \\
\hline
\!\!\!\!-120 & \!\!\! -274,\! -225,\! -85,\! -15\!\!
\end{array}\right]
\end{aligned}
\end{equation*}
$N_{1i} = [0_2^{\top},  1]^{\top}$ and
$N_{2i} = [0_4^{\top}, 1]^{\top}$, we have
$\Psi_{1i} = [3, 6, 5]$ and
$\Psi_{2i} = [120, 270, 225, 80, 15]$.
Finally, we derive the distributed controller (\ref{cotroller}).

Fig. 2 (a) shows the trajectory of $\log(\Vert \mathbf{p}(t) - \mathbf{p}^*\Vert)$ under (\ref{gradient-play}), and indicates  $\mathbf{p}(t)$ converges to $\mathbf{p}^*$ with an exponential rate.
Fig. 2(b) presents the trajectories of $e_i(t), i \in \{1, \dots, 4\}$, where $e_i(t) = y_i(t) - p_i(t)$, and implies that Problem \ref{prb} can be solved by (\ref{cotroller}) since $\lim_{t \to \infty} e_i(t) = 0$.
\begin{figure}[htp]
\centering
\includegraphics[scale=0.3]{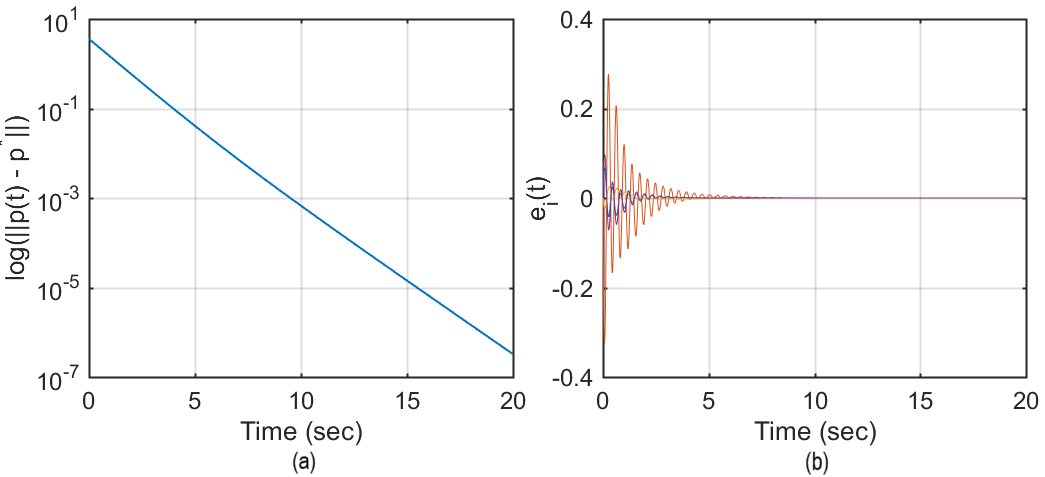}
\caption{(a) The trajectory of $\log(\Vert \mathbf{p}(t) - \mathbf{p}^*\Vert)$.
	(b) The trajectories of $e_i(t)$.}
\label{fig:sim}
\end{figure}

\section{CONCLUSION}
\label{sec:conclusion}

This paper investigated seeking an NE of a monotone game over a  multi-agent system with each agent represented by a nonlinear uncertain dynamics in a lower-triangular form.
Resorting to a reference signal generator to find an NE, and internal models to handle external disturbances, the problem was cast into a robust stabilization problem of an augmented system. Under a set of standard assumptions, the augmented system was semi-globally stabilized by a linear distributed state-feedback controller, which led to the solution of our problem.
Numerical simulations were carried out for illustration.

\bibliographystyle{IEEEtran}
\bibliography{references.bib}

\end{document}